\documentclass{amsart}
\usepackage{amssymb,latexsym,graphics,amsfonts}

\swapnumbers
 \newtheorem{theo}{Theorem}[section]
\newtheorem{cor}[theo]{Corollary}
\newtheorem{lemm}[theo]{Lemma}
\newtheorem{prop}[theo]{Proposition}
\theoremstyle{definition}
\newtheorem{defi}[theo]{Definition}
\theoremstyle{definition}
\newtheorem{ex}[theo]{Example}

\numberwithin{equation}{section}
 \numberwithin{equation}{subsection}

\begin{document}

\title{Module Biprojective and Module Biflat Banach Algebras}

\author[A. Bodaghi]{Abasalt Bodaghi}
\address{Department of Mathematics, Islamic Azad University, Garmsar Branch, Garmsar,
 Iran}
\email{abodaghi@iau-garmsar.ac.ir}

\author[M. Amini]{Massoud Amini}
\address{Department of Mathematics,
 Tarbiat Modares University, Tehran 14115-134, Iran}
\email{mamini@modares.ac.ir}



\dedicatory{dedicated to Professor Alireza Medghalchi}
\smallskip

\begin{abstract}

In this paper we define module biprojctivity and module biflatness
for a Banach algebra which is a Banach module over another Banach
algebra with compatible actions, and find their relation to
classical biprojectivity and biflatness. As a typical example, We
show that for an inverse semigroup $S$ with an upward directed set
of idempotents $E$, the semigroup algebra $ \ell ^{1}(S)$, as an
$\ell ^{1}(E)$-module, is module biprojective if and only if an
appropriate group homomorphic image of $S$ is finite. Also we show
that $\ell ^{1}(S)$ is module biflat if and only if $S$ is amenable.
\end{abstract}

\maketitle

\section{Introduction}
For a discrete semigroup $S$, $\ell ^{\infty}(S)$ is the Banach
algebra of bounded complex-valued functions on $S$ with the supremum
norm and pointwise multiplication. For each $a \in S$ and $f \in
\ell ^{\infty}(S)$, let $l_af$ and $r_af$ denote the left and the
right translations  of $f$ by $a$, that is $(l_af)(s)=f(as)$ and
$(r_af)(s)=f(sa)$, for each $s \in S$. Then a linear functional
$m\in (\ell ^{\infty}(S))^*$ is called a {\it mean} if $\|
m\|=\langle m,1\rangle=1$;  $m$ is called a {\it left} ({\it right})
{\it invariant mean} if $m(l_af)=m(f) \,(m(r_af)=m(f)$,
respectively) for all $s \in S$ and  $f \in \ell ^{\infty}(S)$. A
discrete semigroup $S$ is called {\it amenable} if there exists a
mean $m$ on $\ell ^{\infty}(S)$ which is both left and right
invariant (see \cite{dun}). An {\it inverse semigroup} is a discrete
semigroup $S$ such that for each $s\in S$, there is a unique element
$s^*\in S$ with $ss^*s=s$ and $s^*ss^*=s^*$. Elements of the form
$ss^*$ are called {\it idempotents} of $S$. For an inverse semigroup
$S$, a left invariant mean on $\ell ^{\infty}(S)$ is right invariant
and vise versa.

A Banach algebra $\mathcal A$ is {\it amenable} if every bounded
derivation from $\mathcal A$ into any dual Banach $A$-module is
inner, equivalently if $H^1(\mathcal A,X^*)=\{0\}$ for every
Banach $A$-module $X$, where $H^1(\mathcal A,X^*)$ is the {\it
first Hochschild cohomology group} of $A$ with coefficients in
$X^*$. This concept was introduced by Barry Johnson in
\cite{joh}. Also ${\mathcal A}$ is called {\it super-amenable}
({\it contractible}) if $H^1({\mathcal A},X)=\{0\}$ for every
Banach ${\mathcal A}$-bimodule $X$ (see \cite{dal, run}).

The second author in \cite{am1} introduced the concept of module
amenability for Banach algebras which are Banach modules on another
Banach algebra with compatible actions, and showed that for an
inverse semigroup $S$ with set of idempotents $E$, the semigroup
algebra $\ell^1(S)$ is module amenable as a Banach module on
$\ell^1(E)$, if and only if $S$ is amenable. This generalizes the
celebrated Johnson's theorem for a discrete group $G$ (valid also
for locally compact groups) which states that the group algebra
$\ell^1(G)$ is amenable if and only if $G$ is amenable. The authors
in \cite{boa} introduced the concept of module super-amenability and
showed that for an inverse semigroup $S$, the semigroup algebra
$\ell^1(S)$ is module super-amenable if and only if the group
homomorphic image $S/\approx$ is finite, where $s\approx t$ whenever
$\delta_s-\delta_t$ belongs to the closed linear span of the set
$$\{\delta_{set}-\delta_{st}: s,t\in S,e\in E\}.$$

Biprojective Banach algebras were introduce by A. Ya. Helemskii in
\cite{hel2}. Later he has studied biprojectivity and biflatness of
the Banach algebras in more details in \cite[Chapters IV and
VII]{hel}. It follows from \cite[Proposition 2.8.41]{dal} that a
biprojective Banach algebra is biflat, but the converse is not true.
For an infinite-dimensional Hilbert space $\mathcal H$, the Banach
algebra $A(\mathcal H \widehat \bigotimes \mathcal H)$ consisting of
norm limits of all sequences of finite rank operators on $\mathcal H
\widehat \bigotimes \mathcal H$ is biflat, but not biprojective
\cite[Example 4.3.25]{run}. If $\mathcal K$ is a non-empty, locally
compact space, then $C_0(\mathcal K)$ is biprojective if and only if
$\mathcal K$ is discrete \cite[Proposition 4.2.31]{dal}, but if
$\mathcal K$ is an infinite compact space, then $C(\mathcal K)$ is
not biprojective \cite[Corollary 5.6.3]{dal}. In general each
commutative $C^*$-algebra with discrete character space is
biprojective and the converse holds for all commutative Banach
algebras \cite[Exercise 4.3.5]{run}. It is shown by Selivanov in
\cite{sel}, that for any a non-zero Banach space $E$, the nuclear
algebra $ E \widehat \bigotimes E^*$ is biprojective (see
\cite[Corollary 2.8.43]{dal}).

In part two of this paper, we define the module biprojctivity and
module biflatness of a Banach algebra $\mathcal A$ which is a Banach
$\mathfrak A$-module with compatible actions on another Banach
algebra $\mathfrak A$ and find their relation with module
amenability \cite{am1} and module super-amenability \cite{boa}. In
particular, we show that when $\mathfrak A$ acts on $\mathcal A$
trivially from left then under some mild conditions, module
biprojectivity (biflatness) of $\mathcal A$ implies biprojectivity
(biflatness) of the quotient Banach algebra $\mathcal A/J$ , where
$J$ is the closed ideal of $\mathcal A$ generated by
$\alpha\cdot(ab)-(ab)\cdot\alpha$ for all $a\in \mathcal A$ and
$\alpha\in \mathfrak A$. Also we show that, under some conditions,
biprojectivity (biflatness) of a Banach algebra implies the module
biprojectivity (biflatness) of it, but the converse is not true.

It is shown in \cite[Proposition 3.6]{abe} that if $\mathcal A$ is a
commutative Banach ${\mathfrak A}$-module such that $\mathcal
A^{**}$ is ${\mathfrak A}$-module amenable, then $\mathcal A$ is
module amenable. In part three of this paper we improve this result
by assuming a weaker condition on $\mathcal A$, using our results on
module biflatness.

Let $P$ be a partially ordered set. For $p\in P$, we set
$(p]=\{x:x\leq p\}$ and $[p)=\{x:p\leq x\}$. Then $P$ is called {\it
locally finite} if $(p]$ is finite for all $p\in P$, and {\it
locally $C$-finite} for some constant $C>1$ if $|(p]|<C$ for all
$p\in P$. A partially ordered set $P$ which is locally $C$-finite,
for some constant $C$ is called {\it uniformly locally finite}. Y.
Choi in \cite[Theorem 6.1]{cho} proved that if $S$ is a Clifford
semigroup, then the group algebra $\ell^1(S)$ is biflat if and only
if $(E,\leq)$ is uniformly locally finite, and each maximal subgroup
of $S$ is amenable. Later, P. Ramsden generalized this result to any
discrete semigroup $S$ \cite{ram}. He also showed that for any
discrete semigroup $S$,  $\ell^1(S)$ is biprojective if and only if
$S$ is uniformly locally finite and and all maximal subgroup of $S$
are finite.

In part four of this paper, we prove that if $S$ is an inverse
semigroup with an
 upward directed set of idempotents $E$, then $\ell ^{1}(S)$ is
module biprojective, as an $\ell ^{1}(E)$-module, if and only if an
appropriate group homomorphic image ${S}/{\approx}$ of $S$ is
finite. This could be considered as the module version (for inverse
semigroups) of a result of Helemskii \cite{hel} which asserts that
for a discrete group $G$, $\ell^1(G)$ is biprojective if and only if
$G$ is finite (see also \cite[Theorem 3.3.32]{dal}). Finally we show
that $\ell ^{1}(S)$ is $\ell ^{1}(E)$-module biflat if and only if
$S$ is amenable.

\section{module Biprojective and module Biflat  of Banach algebras}

Throughout this paper, ${\mathcal A}$ and ${\mathfrak A}$ are
Banach algebras such that ${\mathcal A}$ is a Banach ${\mathfrak
A}$-bimodule with compatible actions, that is

$$
\alpha\cdot(ab)=(\alpha\cdot a)b,
\,\,(ab)\cdot\alpha=a(b\cdot\alpha) \hspace{0.3cm}(a,b \in
{\mathcal A},\alpha\in {\mathfrak A}).
$$
Let ${\mathcal X}$ be a Banach ${\mathcal A}$-bimodule and a
 Banach ${\mathfrak A}$-bimodule with compatible actions, that is
$$
\alpha\cdot(a\cdot x)=(\alpha\cdot a)\cdot x,
\,\,a\cdot(\alpha\cdot x)=(a\cdot\alpha)\cdot x, \,\,(\alpha\cdot
x)\cdot a=\alpha\cdot(x\cdot a) \hspace{0.3cm}(a \in{\mathcal
A},\alpha\in {\mathfrak A},x\in{\mathcal X} )
$$
and the same for the right or two-sided actions. Then we say that
${\mathcal X}$ is a Banach ${\mathcal A}$-${\mathfrak A}$-module.
If moreover
$$\alpha\cdot x=x\cdot\alpha \hspace{0.3cm}( \alpha\in {\mathfrak
A},x\in{\mathcal X} )$$ then $\mathcal X $ is called a {\it
commutative} ${\mathcal A}$-${\mathfrak A}$-module. If $\mathcal X
$ is a (commutative) Banach ${\mathcal A}$-${\mathfrak
A}$-module, then so is $\mathcal X^*$, where the actions of
$\mathcal A$ and ${\mathfrak A}$ on $\mathcal X^*$ are defined by
$$\langle\alpha\cdot f,x\rangle{}=\langle{}f,x\cdot\alpha\rangle{},\,\,\langle{}
a\cdot f,x\rangle{}=\langle{}f,x\cdot a\rangle{}\hspace{0.3cm} (a
\in{\mathcal A},\alpha\in {\mathfrak A},x\in{\mathcal X},f \in
\mathcal X^* )$$
 and the same for the right actions. Let ${\mathcal Y}$ be another ${\mathcal A}$-${\mathfrak
A}$-module, then a ${\mathcal A}$-${\mathfrak A}$-module morphism
from ${\mathcal X}$ to ${\mathcal Y}$ is a norm-continuous map
$\varphi :{\mathcal X}\longrightarrow {\mathcal Y}$ with $
\varphi(x \pm y)= \varphi (x)\pm
 \varphi (y)$ and
 $$ \varphi (\alpha \cdot x)= \alpha \cdot\varphi(x), \,\,\varphi(x\cdot \alpha)= \varphi (x)\cdot\alpha
,\,\,\varphi (a\cdot x)=a\cdot \varphi(x), \varphi(x\cdot a)=
\varphi(x)\cdot a,$$ for $x,y \in {\mathcal X}, a \in {\mathcal
A},$ and $\alpha \in {\mathfrak A}$.

 Note that when
${\mathcal A}$ acts on itself by algebra multiplication, it is
not in general a Banach ${\mathcal A}$-${\mathfrak A}$-module, as
we have not assumed the compatibility condition
$$a\cdot(\alpha\cdot b)=(a\cdot\alpha)\cdot b\quad (\alpha\in {\mathfrak A}, a,b \in{\mathcal
A}).$$ If $\mathcal A$ is a commutative $\mathfrak A$-module and
acts on itself by multiplication from both sides, then it is also
a Banach ${\mathcal A}$-${\mathfrak A}$-module.

If ${\mathcal A}$ is a Banach $\mathfrak A$-module with
compatible actions, then so are the dual space $\mathcal A^*$ and
the second dual space $\mathcal A^{**}$. If moreover $\mathcal A$
is a commutative $\mathfrak A$-module, then $\mathcal A^*$ and the
$\mathcal A^{**}$ are commutative ${\mathcal A}$-${\mathfrak
A}$-modules. Also the canonical embedding $\hat{}: \mathcal A\to
\mathcal A^{**}; a\mapsto \hat a$ is an $\mathfrak A$-module
morphism.

Consider the projective tensor product $\mathcal A \widehat
\bigotimes \mathcal A$. It is well known that $\mathcal A
\widehat \bigotimes \mathcal A$ is a Banach algebra with respect
to the canonical multiplication map defined by
$$(a\otimes b)(c\otimes d)=(ac\otimes bd)$$
and extended by bi-linearity and continuity \cite{dal}. Then
$\mathcal A \widehat \bigotimes \mathcal A$ is a Banach ${\mathcal
A}$-${\mathfrak A}$-module with canonical actions. Let $I$ be the
closed ideal of the projective tensor product $\mathcal A \widehat
\bigotimes \mathcal A$ generated by elements of the form $\alpha
\cdot a \otimes b-a \otimes b\cdot\alpha$ for $ \alpha\in {\mathfrak
A},a,b\in{\mathcal A}$. Consider the map $\omega \in {\mathcal
L}(\mathcal A \widehat \bigotimes \mathcal A, \mathcal A )$ defined
by $\omega (a \otimes b)=ab$ and extended by linearity and
continuity. Let $J$ be the closed ideal of ${\mathcal A}$ generated
by $\omega(I)$. Then the module projective tensor product ${\mathcal
A}\widehat \bigotimes _{\mathfrak A} {\mathcal A}\cong(\mathcal A
\widehat \bigotimes \mathcal A)/{I}$ and the quotient Banach algebra
$\mathcal A/J$ are Banach ${\mathfrak A}$-modules with compatible
actions. We have $({\mathcal A}\widehat \bigotimes _{\mathfrak A}
{\mathcal A})^*={\mathcal L}_{\mathfrak A}({\mathcal A},{\mathcal
A^*})$ where the right hand side is the space of all ${\mathfrak
A}$-module morphism from $A$ to $A^*$\cite{rie}. Also the map
$\widetilde{\omega} \in {\mathcal L}({\mathcal A}\widehat \bigotimes
_{\mathfrak A} {\mathcal A},\mathcal A/J)$ defined by $
\widetilde{\omega} (a \otimes b +I)=ab+J$ extends to an ${\mathfrak
A}$-module morphism. Moreover, $ \widetilde{\omega}^*$
 and $\widetilde{\omega}^{**}$, the first and second adjoints of $ \widetilde{\omega} $
  are ${\mathcal A}$-${\mathfrak A}$-module morphisms. If $\mathcal
  A/J$ and ${\mathcal A}\widehat \bigotimes _{\mathfrak A} {\mathcal
  A}$ are commutative $\mathfrak A$-module, then $\mathcal
  A/J$ and ${\mathcal A}\widehat \bigotimes _{\mathfrak A} {\mathcal
  A}$ are ${\mathcal A}/J$-${\mathfrak A}$-module and $ \widetilde{\omega} $ and $\widetilde{\omega}^{**}$ are
${\mathcal A}/J$-${\mathfrak A}$-module homomorphisms. Let $
\square$ and $ \lozenge$ be the first and second Arens products
on the second dual space $\mathcal A^{**}$, then $\mathcal
A^{**}$ is a Banach algebra with respect to both of these
products \cite[Theorem 2.6.15]{dal}.

Let ${\mathcal A}$ and ${\mathfrak A}$ be as in the above  and
${\mathcal X}$ be a Banach ${\mathcal A}$-${\mathfrak A}$-module.
Let $I$ and $J$ be the corresponding closed ideals of $\mathcal A
\widehat \bigotimes \mathcal A$ and $\mathcal A$, respectively. A
bounded map $D: \mathcal A \longrightarrow \mathcal X $ is called
a {\it module derivation} if
$$D(a\pm b)=D(a)\pm D(b),\hspace{0.2cm}D(ab)=D(a)\cdot b+a\cdot D(b)\hspace{0.3cm}
(a,b \in \mathcal A),$$and
$$D(\alpha\cdot a)=\alpha\cdot D(a),\hspace{0.3cm}D(a\cdot\alpha)=D(a)\cdot\alpha
\hspace{0.3cm}(a \in{\mathcal A},\alpha\in {\mathfrak A}).$$
Although $D$ is not necessary linear, but still its boundedness
implies its norm continuity (since it preserves subtraction).
When $\mathcal X $ is commutative, each $x \in \mathcal X $
defines a module derivation
$$D_x(a)=a\cdot x-x\cdot a \hspace{0.3cm} (a \in{\mathcal A}).$$
These are called {\it inner} module derivations. The Banach
algebra ${\mathcal A}$ is called {\it module amenable} (as an
${\mathfrak A}$-module) if for any commutative Banach ${\mathcal
A}$-${\mathfrak A}$-module $\mathcal X $, each module derivation
$D: \mathcal A \longrightarrow \mathcal X^*$ is inner \cite{am1}.
Similarly, ${\mathcal A}$ is called {\it module super-amenable}
if each module derivation $D: \mathcal A \longrightarrow \mathcal
X$ is inner \cite{boa}.

Recall that a Banach algebra $\mathcal A$ is called biprojective
if ${\omega} $ has a bounded right inverse which is an
 ${\mathcal A}$-bimodule homomorphism, and is called biflat if ${\omega}^*$ has a bounded left inverse
which is an ${\mathcal A}$-bimodule homomorphism.

\begin{defi}\label{def1} A Banach algebra $\mathcal A$ is called {\it module
biprojective} (as an $\mathfrak A$-module) if $ \widetilde{\omega} $
has a bounded right inverse which is an
 ${\mathcal A}/J$-${\mathfrak A}$-module morphism.
\end{defi}
\begin{defi}\label{def2} A Banach algebra $\mathcal A$ is called {\it module biflat}
(as an $\mathfrak A$-module) if $ \widetilde{\omega}^* $ has a
bounded left inverse which is an ${\mathcal A}/J$-${\mathfrak
A}$-module morphism.
\end{defi}

We say the Banach algebra ${\mathfrak A}$ acts trivially on
$\mathcal A$ from left if for each $\alpha\in \mathfrak A$ and $a\in
\mathcal A$, $\alpha\cdot a=f(\alpha)a$, where $f$ is a continuous
linear functional on ${\mathfrak A}$. The following lemma is proved
in \cite[Lemma 3.1]{abe}.

\begin{lemm} \label{tpcl} {\it Let ${\mathfrak
A}$ acts on $\mathcal A$ trivially from left (right) and $J_0$ be
a closed ideal of $\mathcal A$ such that $J \subseteq J_0$. If
$\mathcal A/J_0$ has a left (right) identity $e+J_0$, then  for
each $\alpha \in {\mathfrak A}$ and $a\in\mathcal A$
 we have $a\cdot\alpha-f(\alpha)a \in J_0 \hspace{0.2cm}(\alpha\cdot a-f(\alpha)a \in J_0, respectively)$}.
\end{lemm}

 By the definition of
$I$, if $x=\sum_{i=1}^na_i\otimes b_i\in \mathcal A \widehat
\bigotimes \mathcal A$, then $\alpha\cdot x-x\cdot\alpha \in I$ for
all $\alpha\in \mathfrak A$. Therefore $\alpha\cdot
x+I=x\cdot\alpha+I$, that is ${\mathcal A}\widehat \bigotimes
_{\mathfrak A} {\mathcal A}$ is always a commutative $\mathfrak
A$-module. We use this fact in the proof of the following
proposition.

\begin{prop} \label{pro}{\it Assume that ${\mathfrak
A}$ acts on $\mathcal A$ trivially from left and $\mathcal A/J$ has
a left identity. If $\mathcal A$ is biprojective, then $\mathcal A$
is module biprojective.} \end{prop}

\vspace{.3cm}\paragraph{\large\bf Proof.} suppose that $\rho$ is the
the bounded right inverse of ${\omega} $. Define
$\widetilde{\rho}:\mathcal A/J\longrightarrow (\mathcal A \widehat
\bigotimes \mathcal A)/{I}$ via $\widetilde{\rho}(a +J):=
\rho(a)+I$, for all $a\in \mathcal A$. The map $\rho$ is a $\mathcal
A$-bimodule homomorphism, so for each $a,b\in \mathcal A$ and
$\alpha\in \mathfrak A$ we have
\begin{align*}
\rho(\alpha\cdot ab-ab\cdot\alpha) &=\rho((\alpha\cdot
a)b-a(b\cdot\alpha)) \\
&=(\alpha\cdot a)\rho(b)-\rho(a)(b\cdot\alpha)\\
&=\alpha\cdot (a\cdot\rho(b))-(\rho(a)\cdot b)\cdot\alpha\\
&=\alpha\cdot \rho(ab)-\rho(ab)\cdot\alpha.
\end{align*}

Hence $\rho(\alpha\cdot ab-ab\cdot\alpha)\in I$. Thus
$\widetilde{\rho}$ is well defined. If
$\rho(a)=\sum_{i=1}^nx_i\otimes y_i$, where $x_i,y_i\in \mathcal
A$, then
\begin{align*}
\rho(\alpha\cdot a)&=\rho(f(\alpha)a)
=f(\alpha)\rho(a)\\
&=f(\alpha)\sum_{i=1}^nx_i\otimes
y_i=\sum_{i=1}^nf(\alpha)x_i\otimes y_i\\
&=\sum_{i=1}^n\alpha\cdot x_i\otimes y_i=\alpha
\cdot\sum_{i=1}^nx_i\otimes y_i=\alpha\cdot\rho(a).
\end{align*}
Hence $\widetilde{\rho}(\alpha\cdot(
a+J))=\alpha\cdot\widetilde{\rho}(a+J)$. Also by Lemma \ref{tpcl} we
have
\begin{align*}
\widetilde{\rho}(( a+J)\cdot\alpha)&=\widetilde{\rho}(
a\cdot\alpha+J)=\widetilde{\rho}(f(\alpha)a+J)\\
&=f(\alpha)\rho(a)+I=\alpha\cdot\rho(a)+I\\
&=\rho(a)\cdot\alpha+I=\widetilde{\rho}(a+J)\cdot\alpha.
\end{align*}

Obviously $\widetilde{\rho}$ is a $\mathcal A/J$-bimodule
homomorphism. Hence $\widetilde{\rho}$ is a ${\mathcal
A}/J$-${\mathfrak A}$-module morphism. Now for each $a\in \mathcal
A$ we have $(
\widetilde{\omega}\circ\widetilde{\rho})(a+J)=\widetilde{\omega}(\rho(a)+I)=\omega(\rho(a))+J=a+J$.
Therefore $\widetilde{\rho}$ is a bounded right inverse for
$\widetilde{\omega}$.$\hfill\blacksquare $

\begin{prop} \label{bflat} {\it Assume that ${\mathfrak
A}$ acts on $\mathcal A$ trivially from left and $\mathcal A/J$ has
a left identity. If $\mathcal A$ is biflat, then $\mathcal A$ is
module biflat.} \end{prop}

 \vspace{.3cm}\paragraph{\large\bf Proof.} It is straightforward to show that $\widetilde{\omega}\circ\pi_1=\pi_2\circ\omega$,
where $\pi_1:\mathcal A \widehat \bigotimes \mathcal
A\longrightarrow (\mathcal A \widehat \bigotimes \mathcal A)/I$
and $\pi_2:\mathcal A\longrightarrow \mathcal A/J$ are projection
maps. Assume that $\theta$ is the the bounded left inverse of
${\omega}^* $.

Define $\widetilde{\theta}:(\mathcal A \widehat \bigotimes
\mathcal A)/{I})^*\longrightarrow (\mathcal A/J)^*$ via
$$(\widetilde{\theta}(\phi))(a+J):= [\theta(\phi\circ\pi_1)](a)\quad (a\in \mathcal A),$$
where $\phi$ is a functional in $((\mathcal A \widehat \bigotimes
\mathcal A)/{I})^*$. As in the proof of Proposition \ref{pro}, we
can show that $\widetilde{\theta}$ is a well defined  ${\mathcal
A}/J$-${\mathfrak A}$-module morphism. If $f$ is a bounded
functional on $\mathcal A/I$ and $a\in \mathcal A$, then
\begin{align*}
[(\widetilde{\theta}\circ \widetilde{\omega}^*)(f)](a+J)&=(\widetilde{\theta}( \widetilde{\omega}^*)(f))(a+J) \\
&= [(\widetilde{\theta}( \widetilde{\omega}^*(f)\circ\pi_1)](a) \\
&= [(\theta(f\circ \widetilde{\omega}\circ\pi_1)](a)\\
&= [(\theta(f\circ\pi_2\circ\omega)](a)\\
&= [(\theta({\omega}^*(f\circ\pi_2))](a)\\
&= (f\circ\pi_2)(a)=f(a+J).
\end{align*}
Therefore $\widetilde{\theta}$ is a bounded right inverse for
$\widetilde{\omega}^*$.$\hfill\blacksquare$

In section 4 we give an example of a Banach algebra which is module
biprojective (biflat), but not biprojective (biflat). We show that
module biprojectivity (biflatness) of $\mathcal A$ implies
biprojectivity (biflatness) of $\mathcal A/J$.

\begin{prop} \label{prom}{\it Let ${\mathfrak
A}$ acts on $\mathcal A$ trivially from left and $\mathcal A/J$
has a left identity such that ${\mathfrak A}$ has a bounded
approximate identity for $\mathcal A$. If $\mathcal A$ is module
biprojective, then $\mathcal A/J$ is biprojective.} \end{prop}

\vspace{.3cm}\paragraph{\large\bf Proof.} Suppose that
$\widetilde{\rho}$ is the the bounded right inverse of
$\widetilde{\omega} $. We show that the  map
$\overline{\omega}:\mathcal A/J\widehat\bigotimes\mathcal
A/J\longrightarrow\mathcal A/J;\quad
((a+J)\otimes(b+J)\longrightarrow ab+J)$ has a right inverse.
Consider the map $\Gamma:(\mathcal A\widehat\bigotimes\mathcal
A)/ker (\pi\otimes\pi)\longrightarrow \mathcal
A/J\widehat\bigotimes\mathcal A/J;\quad (a \otimes b+ker
(\pi\otimes\pi)\longrightarrow (a+J)\otimes (b+J))$, where
$\pi:\mathcal A\longrightarrow \mathcal A/J$ is the projection
map. We have $I\subseteq ker (\pi\otimes\pi$) because for each
$a,b\in \mathcal A$ and $\alpha\in \mathfrak A$
\begin{align*}
(\pi\otimes\pi)(\alpha \cdot a \otimes b-a \otimes b\cdot\alpha)&=
(\alpha \cdot a+J) \otimes (b+J)-(a+J) \otimes (b\cdot\alpha+J)\\
&= (f(\alpha)a+J) \otimes (b+J)-(a+J) \otimes (f(\alpha)b+J)\\
&= f(\alpha)(a+J) \otimes (b+J)-f(\alpha)(a+J) \otimes (b+J)=0.
\end{align*}
We have used Lemma \ref{tpcl} in the second equality. Hence the map
 $$\Phi:(\mathcal
A\widehat\bigotimes\mathcal A)/I\longrightarrow (\mathcal
A\widehat\bigotimes\mathcal A)/ker (\pi\otimes\pi)\quad
(x+I\longrightarrow x+ ker (\pi\otimes\pi))$$
 is well defined. We put
$\overline{\rho}=\Gamma\circ\Phi\circ\widetilde{\rho}$. Since
${\mathfrak A}$ has a bounded approximate identity for $\mathcal A$,
wit follows that $\overline{\rho}$ is $\mathbb{C}$-linear. Now if
$\widetilde{\rho}(a+J)=\sum_{i=1}^nx_i\otimes y_i+I$, where
$x_i,y_i\in \mathcal A$, then
\begin{align*}
\langle\overline{\omega}\circ\overline{\rho}, a+J\rangle&=
\langle\overline{\omega}\circ
\Gamma\circ\Phi\circ\widetilde{\rho}, a+J\rangle\\
&= \langle\overline{\omega}\circ
\Gamma\circ\Phi, \sum_{i=1}^nx_i\otimes y_i+I\rangle\\
&= \langle\overline{\omega}\circ \Gamma, \sum_{i=1}^nx_i\otimes y_i+ ker (\pi\otimes\pi)\rangle\\
&= \langle\overline{\omega}, \sum_{i=1}^n(x_i+J)\otimes
(y_i+J)\rangle.\\
&= \sum_{i=1}^n\omega(x_i\otimes y_i)+J\\
&= \widetilde{\omega}(\sum_{i=1}^n(x_i\otimes y_i)+I)\\
&= \widetilde{\omega}\circ\widetilde{\rho}(a+J)=a+J.
\end{align*}
Therefore $\overline{\rho}$ is right inverse of
$\overline{\omega}$.$\hfill\blacksquare$

We say that ${\mathfrak A}$ has a bounded approximate identity for
$\mathcal A$ if there is a bounded net $\{\alpha_j \}$ in $\mathfrak
A$ such that $\|\alpha_j.a-a\|\to 0$ and $\|a.\alpha_j-a\|\to 0$,
for each  $a\in {\mathcal A}$.

\begin{prop} \label{bifm}{\it Assume that ${\mathfrak
A}$ acts on $\mathcal A$ trivially from left, $\mathcal A/J$ has a
left identity, and ${\mathfrak A}$ has a bounded approximate
identity for $\mathcal A$. If $\mathcal A$ is module biflat, then
$\mathcal A/J$ is biflat.} \end{prop}

\vspace{.3cm}\paragraph{\large\bf Proof.} Assume that $\Phi$,
$\Gamma$ and $\overline{\omega}$ are as the above. Suppose that
$\widehat{\rho}$ is the the bounded right inverse of
$\widetilde{\omega}^*$. We prove that the map
$\overline{\omega}^*:(\mathcal A/J)^*\longrightarrow (\mathcal
A/J\widehat\bigotimes\mathcal A/J)^*$ has a left inverse. From the
proof of Proposition \ref{prom} we see that $\overline{\omega}\circ
\Gamma\circ\Phi=\widetilde{\omega}$. Now for each $\varphi\in
(\mathcal A/J)^*$ we have
\begin{align*}
(\widehat{\rho}\circ(\Gamma\circ\Phi)^*\circ\overline{\omega}^*)(\varphi)&=
(\widehat{\rho}\circ(\Gamma\circ\Phi)^*)(\varphi\circ\overline{\omega})\\
&= \widehat{\rho}(\varphi\circ\overline{\omega}\circ
\Gamma\circ\Phi)\\
&= (\widehat{\rho}\circ\widetilde{\omega}^*)(\varphi)=\varphi.
\end{align*}
Therefore the map $\widehat{\rho}\circ(\Gamma\circ\Phi)^*$ is left
inverse of $\overline{\omega}^*$ which is $\mathbb{C}$-linear
(see the proof of Proposition \ref{prom}). $\hfill\blacksquare$

One can easily show that module biprojectivity implies module
biflatness. In section 4, we shall give an example of a Banach
algebra which is module biflat but not module biprojective.

Let $X$, $Y$ and $Z$ be Banach ${\mathcal A}/J$-${\mathfrak
A}$-modules. Then the short, exact sequence

$\hspace{5cm}\{0\}\longrightarrow X
\stackrel{\varphi}{\longrightarrow}Y
\stackrel{\psi}{\longrightarrow} Z\longrightarrow \{0\}\hfill (2.1)$
\linebreak is {\it admissible} if $\psi$ has a bounded right inverse
which is ${\mathfrak A}$-module homomorphism, and {\it splits} if
$\psi$ has a bounded right inverse which is a ${\mathcal
A}/J$-${\mathfrak A}$-module morphism. Obviously, the short, exact
sequence (2.1) is admissible if and only if $\varphi$ has a bounded
left inverse which is ${\mathcal A}/J$-${\mathfrak A}$-module
morphism.

We set $K=\ker(\widetilde{\omega})$. If ${\mathcal A}/J$ has a
bounded approximate identity, then the following sequences are
exact.

$\hspace{4cm}\{0\}\longrightarrow K
\stackrel{i}{\longrightarrow}(\mathcal A \widehat \bigotimes
\mathcal A)/{I} \stackrel{\widetilde{\omega}}{\longrightarrow}
{\mathcal A}/J\longrightarrow \{0\}\hfill (2.2)$\\

$\hspace{3.5cm}\{0\}\longrightarrow ({\mathcal A}/J)^*
\stackrel{\widetilde{\omega}^*}{\longrightarrow}
 (\mathcal A \widehat \bigotimes
\mathcal A)/{I})^*\stackrel{i^*}{\longrightarrow}
K^*\longrightarrow \{0\}\hfill (2.3)$

\begin{defi} {\it A
bounded net $\{\widetilde{\xi_j} \}$ in ${\mathcal A}\widehat
\bigotimes _{\mathfrak A} {\mathcal A}$ is called a {\it module
approximate diagonal} if $\widetilde{\omega}_{\mathcal A}(
\widetilde{\xi_j}) $ is a bounded approximate identity of  $
{\mathcal A}/{J}$ and
$$ \lim_j \| \xi_j \cdot a- a\cdot \xi_j \| =0 \hspace{0.3cm} (a \in \mathcal A
). $$ An element $ \widetilde{E} \in ({\mathcal A}\widehat
\bigotimes _{\mathfrak A} {\mathcal A})^{**}$ is called a {\it
module virtual diagonal} if
 $$ \widetilde{\omega}^{**}_{\mathcal A} (
\widetilde{E})\cdot a= \widetilde{a},\hspace{0.3cm}
\widetilde{E}\cdot a=a\cdot \widetilde{E} \hspace{0.3cm} (a \in
\mathcal A),
$$ where $ \widetilde{a}=a+ J^{\perp \perp }$.}
\end{defi}

The following results are proved in \cite[Lemma 3.2]{hoe} and
\cite[Theorem 3.3]{hoe}.

\begin{lemm} \label{lbai}{\it Let $\mathcal
A/J$ be commutative Banach ${\mathfrak A}$-modules. If $\mathcal
A/J$ has an identity, the exact sequences {\text (2.2)} and
{\text (2.3)} are admissible. If $\mathcal A/J$ has a bounded
approximate identity, the exact sequence \text {(2.3)} is
admissible.} \end{lemm}

\begin{theo}\label{bai} {\it Let $\mathcal
A/J$ and ${\mathcal A}\widehat \bigotimes_{\mathfrak A} {\mathcal
A}$ are commutative Banach ${\mathfrak A}$-modules. The Banach
algebra ${\mathcal A}$ is module amenable if and only if $\mathcal
A/J$ has a bounded approximate identity, and the exact sequence
{\text (2.3)} splits.} \end{theo}

Note that ${\mathcal A}\widehat \bigotimes_{\mathfrak A} {\mathcal
A}$ is always a commutative ${\mathfrak A}$-module (see the
paragraph after Lemma \ref{tpcl}). Therefore this condition in
Theorem \ref{bai} is redundant.

The second author in \cite[Proposition 2.2]{am1} showed that if
$\mathcal A$ is a commutative Banach $\mathfrak A$-module which is
module amenable then it has a bounded approximate identity.
Converse is not true for the semigroup algebra even for the
classical case of groups. For the free group $\mathbb{F}_2$ on
two generators, the group algebra $\ell^1(\mathbb{F}_2)$ has a
bounded approximate identity and even an identity, but it is not
amenable \cite[Example 3.3.62]{dal}. The following corollary show
that the converse of this result holds if $\mathcal A$ is module
biflat.

\begin{cor} \label{biflat}{\it Let $\mathcal
A/J$ be a commutative Banach ${\mathfrak A}$-module. Then $\mathcal
A $ is $\mathfrak A$-module amenable if and only if $\mathcal A $ is
module biflat and ${\mathcal A}/{J}$ has a bounded approximate
identity.} \end{cor}

\vspace{.3cm}\paragraph{\large\bf Proof.} This follows immediately
from Theorem \ref{bai}.$\hfill\blacksquare$

Recall that an element $\mathcal M \in {\mathcal A}\widehat
\bigotimes _{\mathfrak A} {\mathcal A}$ is called a {\it module
diagonal} if $\widetilde {\omega}(\mathcal M)$ is an identity of
$\mathcal A/J$ and $a\cdot \mathcal M=\mathcal M\cdot a$, for all
$a\in \mathcal A$. The authors in \cite[Theorem 2.6]{boa} proved
that if $\mathcal A$ is a Banach $\mathcal A$-${\mathfrak
A}$-module, then $\mathcal A$ is module super-amenable if and
only if $\mathcal A$ has a module diagonal. Also it is shown in
\cite[Proposition 2.3]{boa} that if $\mathcal A$ be a commutative
Banach module super-amenable as an ${\mathfrak A}$-module, then
it is unital. The converse of this proposition is true if
$\mathcal A $ is module biprojective.

\begin{theo} \label{bipro}{\it Let $\mathcal A$ be a Banach ${\mathfrak
A}$-module with trivial left action and let $\mathcal A/J$ be a
commutative Banach ${\mathfrak A}$-module. Then $\mathcal A $ is
module super-amenable if and only if $ {\mathcal A}/{J}$ has an
identity and $\mathcal A $ is module biprojective.}
\end{theo}

\vspace{.3cm}\paragraph{\large\bf Proof.} Assume that $\mathcal A
$ is module super-amenable. Then $\mathcal A/J$ is super-amenable
\cite[Lemma 2.7]{boa}, so ${\mathcal A}/{J}$ has identity
\cite[Theorem 2.8.48]{dal}. Also $\mathcal A$ has a module
diagonal $\mathcal M$ \cite[Theorem 2.6]{boa}. Define
$\widetilde{\rho}:\mathcal A/J\longrightarrow (\mathcal A
\widehat \bigotimes \mathcal A)/{I}$ by $\widetilde{\rho}(a +J):=
a\cdot \mathcal M$, for all $a\in \mathcal A$. Since $a\cdot
\mathcal M=\mathcal M\cdot a$, for all $a,b,c,d\in \mathcal A$
and $\alpha\in \mathfrak A$ we have
\begin{align*}
c\otimes d\cdot(\alpha\cdot(ab)-(ab)\cdot\alpha)
&= c\otimes d\cdot(f(\alpha) ab)-c\otimes d\cdot(ab)\cdot\alpha\\
&= f(\alpha)c\otimes dab-c\otimes (dab)\cdot\alpha\\
&= \alpha\cdot c\otimes dab-c\otimes (dab)\cdot\alpha\in I,
\end{align*}
so $\widetilde{\rho}$ is well-defined. Also for all $a,b\in
\mathcal A$ we have
\begin{align*}
\widetilde{\rho}((a+J)\cdot(b+J))&= ab\cdot \mathcal M\\
&= (a+J)\cdot(b\cdot \mathcal M)\\
&= (a+J)\cdot\rho(b+J),
\end{align*}
 and
\begin{align*}
\widetilde{\rho}((b+J)\cdot(a+J))&= ba\cdot \mathcal M= \mathcal M\cdot ba\\
&= (\mathcal M\cdot b).(a+J)\\
&= (b\cdot \mathcal M)\cdot(a+J)\\
&= \widetilde{\rho}(b+J)\cdot (a+J).
\end{align*}
It follows from Lemma \ref{tpcl} that $\widetilde{\rho}$ is a
$\mathcal A/J$- $\mathfrak A$-module morphism. Also
$(\widetilde{\omega}\circ\widetilde{\rho})(a+J)=\widetilde{\omega}(a\cdot
\mathcal M)= (a+J)\cdot \widetilde{\omega}(\mathcal M)=a+J$.
Therefore $\widetilde{\rho}$ is a bounded right inverse for
$\widetilde{\omega}$. Conversely, if $e+J$ is an identity for
${\mathcal A}/{J}$ and $\widetilde{\rho}$ is a bounded right
inverse for $\widetilde{\omega}$ which is ${\mathcal
A}/J$-${\mathfrak A}$-module morphism, it is easy that show
$\widetilde{\rho}(e+J)$ is a module diagonal for $\mathcal A$.
Now the module super-amenability of $\mathcal A$ follows from
\cite[Theorem 2.6]{boa}. $\hfill\blacksquare$

\section{module amenability of the second dual}
In this section we explore conditions under which $\mathcal A$ is
module amenable. Let $\widetilde{\omega}$ be as in the above
section. Let $I$ and $J$ be the corresponding closed ideals of
$\mathcal A \widehat \bigotimes \mathcal A$ and $\mathcal A$,
respectively. Consider the morphism
$$ \widetilde{\omega}^* : ({\mathcal A}/{J})^*=J^{\perp} \longrightarrow
 ({\mathcal A}\widehat
\bigotimes _{\mathfrak A} {\mathcal A})^*={\mathcal L}_{\mathfrak
 A}({\mathcal A},{\mathcal A^*}),$$
$$\langle[\widetilde{\omega}^*(f)](a),b\rangle=f(ab+J), \hspace{0.3cm}f \in ({\mathcal
A}/{J})^*.$$

Let ${\mathcal A}^{**}\widehat \bigotimes _{\mathfrak A} {\mathcal
A}^{**}$ be the projective module tensor product of ${\mathcal
A}^{**}$ and ${\mathcal A}^{**}$, that is ${\mathcal A}^{**}\widehat
\bigotimes _{\mathfrak A} {\mathcal A}^{**}=({\mathcal A^{**}
\widehat \bigotimes \mathcal A^{**}})/{M}$, where $M$ is a closed
ideal generated by elements of the form $\alpha \cdot F \otimes G-F
\otimes G\cdot\alpha$ for $ \alpha\in {\mathfrak A},F,G\in{\mathcal
A}^{**}$. Consider
 $$\widehat{\omega}: {\mathcal A}^{**}\widehat
\bigotimes _{\mathfrak A} {\mathcal A}^{**} \longrightarrow
{\mathcal A^{**}}/{N}$$
$$ \widehat{\omega}(F \otimes G+M)=F \square G+ N,$$
where $ N $ is the closed ideal of ${\mathcal A}^{**}$ generated by
$\widehat{\omega}(M)$. We know that
$$ \widehat{\omega}^*:({\mathcal A^{**}}/N)^* \longrightarrow
 ({\mathcal A}^{**}\widehat
\bigotimes _{\mathfrak A} {\mathcal A}^{**})^*={\mathcal
L}_{\mathfrak A}(\mathcal A^{**},\mathcal A^{***})$$
$$\langle[\widehat{\omega}^*(\phi)]F,G\rangle=\phi(F\square G+ N)\quad(F,G\in{\mathcal A}^{**}) $$

where $\phi\in {\mathcal A^{**}}/N$. Let $T \in ({\mathcal
A}\widehat \bigotimes _{\mathfrak A} {\mathcal A})^*={\mathcal
L}_{\mathfrak A}({\mathcal A},{\mathcal A^*})$ and $T^*$ ,$T^{**}$
are the first and second conjugates of $T$, then $T^{**} \in
({\mathcal A}^{**}\widehat \bigotimes _{\mathfrak A} {\mathcal
A}^{**})^{*}={\mathcal L}_{\mathfrak A}({\mathcal A^{**}},{\mathcal
A}^{***})$ with $\langle T^{**}(F),G\rangle=\lim_j\lim_k\langle
T(a_j),b_k)\rangle$, where $(a_j),(b_k)$ are bounded nets in $
\mathcal A$
 such that  $ \hat a_j \stackrel{w^*}{\longrightarrow}F$
and $ \hat b_k \stackrel{w^*}{\longrightarrow}G$.

\begin{lemm} Let $\mathcal A$ be a Banach
algebra. Then there exist an $\mathcal A$-$\mathfrak A$-module
morphism $\Lambda: {\mathcal L}_{\mathfrak A}({\mathcal A},{\mathcal
A^*})\longrightarrow {\mathcal L}_{\mathfrak A}({\mathcal
A^{**}},{\mathcal A}^{***})$ such that for every T $\in {\mathcal
L}_{\mathfrak A}({\mathcal A},{\mathcal A^*})$, $F,G \in A^{**}$,
and bounded nets $(a_j),(b_k)\subset \mathcal A $ with $ \hat a_j
\stackrel{J^{\perp}}{\longrightarrow}F$ and $ \hat b_k
\stackrel{J^{\perp}}{\longrightarrow}G$, we have
$$\langle\Lambda (T)F,G\rangle= \lim_j\lim_k \langle T(a_j),b_k\rangle.$$
\end{lemm}

 \vspace{.3cm}\paragraph{\large\bf Proof.} It is sufficient to show that $\Lambda $ with
the above definitions is a $\mathcal A$-$\mathfrak A$-module
morphism. It is easy to see that for $f \in J^{\perp} , \alpha\in
{\mathfrak A}$ we have $f\cdot\alpha,\alpha\cdot f \in J^{\perp}$.
If $\alpha \in \mathfrak A$ then $ b_k\cdot\alpha
\stackrel{J^{\perp}}{\longrightarrow}G\cdot\alpha $, hence
\begin{align*}
\langle\Lambda (\alpha \cdot T)F,G\rangle &=\lim_j\lim_k \langle(\alpha\cdot T)(a_j),b_k\rangle \\
&=\lim_j\lim_k \langle\alpha\cdot T(a_j),b_k\rangle\\
&=\lim_j\lim_k \langle T(a_j),b_k\cdot \alpha\rangle,
\end{align*}
and
\begin{align*}
\langle[\alpha\cdot \Lambda (T)]F,G\rangle &=\langle \Lambda (T)F,G\cdot \alpha\rangle \\
&=\lim_j\lim_k \langle T(a_j),b_k\cdot \alpha\rangle.
\end{align*}
$\hfill\blacksquare $

It is shown in the proof of \cite[Theorem 3.4]{abe} that the map
$\lambda : \mathcal A^{**}/N\longrightarrow \mathcal A^{**}/J^{\perp
\perp}; F+N \longrightarrow F+ J^{\perp \perp}$ is a well defined
bounded $\mathcal A$-$\mathfrak A$-module morphism.

\begin{theo}\label{mflat} {\it If  $\mathcal A^{**}$
is module biflat, then so is $\mathcal A$.} \end{theo}

\vspace{.3cm}\paragraph{\large\bf Proof.} Suppose that $
\widehat{\omega}^* $ has a left inverse $\mathcal A/J$-$\mathfrak
A$-module morphism $ \widehat \theta $, then $ \widehat \theta \circ
\widehat{\omega}^* =id_{({\mathcal A^{**}}/N)^*}$. Assume that $
j:({\mathcal A}/{J})^* \longrightarrow {({\mathcal A^{**}}/{J^{\perp
\perp }})^*}$ is the canonical embedding, and $\widetilde{\omega}$,
$\lambda$ and $\Lambda $ are as the above. Consider the map
$i:{\mathcal A}/{J} \longrightarrow {\mathcal A^{**}}/N$;
($a+J\longrightarrow a+N$). Obviously $i$ is well defined. Let us
show that $\Lambda \circ \widetilde{\omega}^* = \widehat{\omega}^*
\circ\lambda^*\circ j $. Take $\varphi \in ({\mathcal A}/{J})^* $,
and $F,G\in{\mathcal A}^{**}, \hat a_l
\stackrel{J^{\perp}}{\longrightarrow}F$ and $ \hat b_k
\stackrel{J^{\perp}}{\longrightarrow}G$, then

\begin{align*}
\langle[(\Lambda \circ \widetilde{\omega}^*)(\varphi
)](F),G\rangle &=\lim_l\lim_k
\langle(\widetilde{\omega}^*(\varphi))(a_l),b_k\rangle \\
&=\lim_l\lim_k \varphi(a_lb_k+J)
\end{align*}
and
\begin{align*}
\langle[(\widehat{\omega}^*\circ\lambda^* \circ j)(\varphi
)](F),G\rangle &=\langle \lambda^*(j(\varphi)),F\square G+ N\rangle\\
&=\langle j(\varphi)\circ\lambda,F\square G+ N\rangle\\
&=\langle j(\varphi),F\square G+ J^{\perp \perp}\rangle\\
&=\lim_l\lim_k\langle a_lb_k+J,\varphi\rangle \\
&=\lim_l\lim_k \varphi(a_lb_k+J).
\end{align*}
Now consider $\Delta=i^* \circ  \widehat \theta \circ \Lambda $,
then it is easy to see that $\Delta$ is a left inverse for
$\widetilde{\omega}^*$.$\hfill\blacksquare $

\begin{lemm} \label{bail}{\it If  $
{\mathcal A^{**}}/{J^{\perp \perp }}= ({\mathcal A}/{J})^{**}$ has a
bounded approximate identity, then so does
 ${\mathcal A}/{J}$.}\end{lemm}

\vspace{.3cm}\paragraph{Proof.} This follows from \cite[Proposition
28.7]{bdu} and \cite[Lemma 1.1]{glw}.$\hfill\blacksquare $

In the next result we give an alternative proof of \cite[Proposition
3.6]{abe} under a weaker condition. In fact if $\mathcal A$ is a
commutative $\mathfrak A$-module, then $\mathcal A^/J$ and $\mathcal
A^{**}/N$ are commutative $\mathfrak A$-modules. Therefore $\mathcal
A^{**}/J^{\perp \perp}$ is also a commutative $\mathfrak A$-module.

\begin{theo}{\it Let $\mathcal A^/J$ and ${\mathcal A}^{**}/N$ be commutative Banach ${\mathfrak A}$-modules.
If $\mathcal A^{**}$ is module amenable, then so is $\mathcal A$.}
\end{theo}

\vspace{.3cm}\paragraph{\large\bf Proof.} Suppose that $\mathcal
A^{**}$ is module amenable, so $\mathcal A^{**}$ is module biflat
and ${\mathcal A}^{**}/N$ has a bounded approximate identity
$\{E_j+N\}$ by Corollary \ref{biflat}. Since $\lambda$ is
surjective, $\{E_j+J^{\perp \perp}\}$ is a bounded approximate
identity for $\mathcal A^{**}/J^{\perp \perp}$. Now the result
follows from Theorem \ref{mflat}, Lemma \ref{bail} and Corollary
\ref{biflat}.$\hfill\blacksquare$

\section{module Biprojectivity and module Biflatness  of semigroup algebras}

In this section we find conditions on a (discrete) inverse
semigroup $S$ such that the semigroup algebra $\ell^1(S)$ is
$\ell^1(E)$-module biprojective and biflat, where $E$ is the set
of idempotents of $S$, acting on $S$ trivially from left and by
multiplication from right. Throughout this section $S$ is an
inverse semigroup with set idempotent $E$, where the order of $E$
is defined by
$$e\leq d \Longleftrightarrow ed=e \hspace{0.3cm}(e,d \in
E).$$ It is easy to show that $E$ is a (commutative) subsemigroup of
$S$ \cite[Theorem V.1.2]{how}. In particular $ \ell ^{1}(E)$ could
be regard as a subalgebra of $ \ell ^{1}(S)$, and thereby $ \ell
^{1}(S)$ is a Banach algebra and a Banach $ \ell ^{1}(E)$-module
with compatible actions \cite{am1}. Here we let $ \ell ^{1}(E)$ act
on $ \ell ^{1}(S)$ by multiplication from right and trivially from
left, that is
$$\delta_e\cdot\delta_s = \delta_s, \,\,\delta_s\cdot\delta_e = \delta_{se} =
\delta_s * \delta_e \hspace{0.3cm}(s \in S,  e \in E).$$

In this case, the ideal $J$ (see section 2) is the closed linear
span of
$$\{\delta_{set}-\delta_{st} \quad s,t \in S,  e \in E
\}.$$
We consider an equivalence relation on $S$ as follows
$$s\approx t \Longleftrightarrow \delta_s-\delta_t \in J \hspace{0.2cm} (s,t \in
S).$$

Recall that $E$ is called {\it upward directed } if for every $e,f
\in E$ there exist $g \in E$ such that $eg=e$ and $fg=f$. This is
precisely the assertion that $S$ satisfies the $D_1$ condition of
Duncan and Namioka \cite{dun}. It is shown in \cite[Theorem
3.2]{ra1} that if $E$ is upward directed, then the quotient
${S}/{\approx}$ is a discrete group. With the notations of
previous section, $ \ell ^{1}(S)/{J}\cong {\ell ^{1}}(S/\approx)$
is a commutative $ \ell ^{1}(E)$-bimodule with the following
actions
$$\delta_e\cdot(\delta_s+J) = \delta_s+J, \,\,(\delta_s+J)\cdot\delta_e = \delta_{se}+J \hspace{0.3cm}(s \in S,  e \in E).$$

It is shown in the paragraph after Lemma \ref{tpcl}, that ${\mathcal
A}\widehat \bigotimes _{\mathfrak A} {\mathcal A}$ is commutative
$\mathfrak A$-module. In particular ${\ell ^{1}(S)}\widehat
\bigotimes _{\ell ^{1}(E)} {\ell ^{1}(S)}\cong{\ell ^{1}(S\times
S)}/{I}$ is a commutative ${\ell ^{1}(E)}$-module, where  $I$ is the
closed linear span of the set of elements of the form
$\delta_{(set,x)} - \delta_{(st,x)}$, where $s,t,x\in S$ and $e\in
E$.

Helemskii in \cite{hel} shows that, for any locally compact group
$G$, $L^1(G)$ is biprojective if and only if $G$ is compact (see
also \cite[Theorem 3.3.32]{dal}). The following Theorem is the
module version of Helemskii's result for inverse semigroups.

\begin{theo}\label{invpro} {\it Let $S$ be an
inverse semigroup with an upward directed set of idempotents $E$.
 Then $\ell ^{1}(S)$ is  module biprojective (as an $\ell
^{1}(E)$-module with trivial left action) if and only if
${S}/{\approx}$ is finite.}\end{theo}

\vspace{.3cm}\paragraph{\large\bf Proof.}  Since $S/\approx$ is a
(discrete) group, ${\ell ^{1}(S)}/{J}\cong {\ell ^{1}}(S/\approx)$
has an identity. By Theorem \ref{bipro}, $\ell ^{1}(S)$ is module
biprojective if and only if $\ell ^{1}(S)$ is module super-amenable.
It follows from \cite[Theorem 3.2]{boa} that $\ell ^{1}(S)$ is
module super-amenable if and only if $S/\approx$ is
finite.$\hfill\blacksquare $

\begin{cor}{\it Let $S$ be an
inverse semigroup with an upward directed set of idempotents $E$. If
$S/\approx$ is infinite, then $\ell ^{1}(S)$ is not
biprojective.}\end{cor}

\vspace{.3cm}\paragraph{\large\bf Proof.} This follows from
Proposition \ref{pro} and Theorem \ref{invpro}.$\hfill\blacksquare $

\begin{theo} \label{invflat}{\it Let $S$ be an
inverse semigroup with an upward directed set of idempotents $E$.
Then $\ell ^{1}(S)$ is module biflat as an $\ell ^{1}(E)$-module
with trivial left action if and only if $S$ is amenable.}\end{theo}

\vspace{.3cm}\paragraph{\large\bf Proof.} Since ${\ell
^{1}(S)}/{J}\cong {\ell ^{1}}(S/\approx)$ has identity, by Corollary
\ref{biflat}, $\ell ^{1}(S)$ is module biflat if and only if $\ell
^{1}(S)$ is module amenable. It follows from \cite[Theorem 3.1]{am1}
that $\ell ^{1}(S)$ is module amenable if and only if $S$ is
amenable.$\hfill\blacksquare $

\begin{cor}{\it Let $S$ be an
inverse semigroup with an upward directed set of idempotents $E$. If
$\ell ^{1}(S)$ is biflat, then $S$ is amenable.}\end{cor}

\vspace{.3cm}\paragraph{\large\bf Proof.} This is a consequence of
Proposition \ref{bflat} and Theorem
\ref{invflat}.$\hfill\blacksquare $

\begin{ex}
(i) Let $\mathbb{N}$ be the commutative semigroup of positive
integers. It is known that $ \ell ^{1}(\mathbb{N})$ with pointwise
multiplication is biprojective, but with convolution it is not
biprojective \cite[Example 4.1.42 ]{dal}. Now consider $(\mathbb{N},
\vee)$ with maximum operation $m\vee n=max(m,n)$, then each element
of $\mathbb{N}$ is an idempotent, hence $\mathbb{N}/\approx$ is the
trivial group with one element. Thus $\ell ^{1}(\mathbb{N})$ is
module biprojective (as an $\ell ^{1}(\mathbb{N})$-module) by
Theorem \ref{invpro}. Since $\mathbb{N}/\approx$ is amenable,
$\mathbb{N}$ is amenable. Therefore $\ell ^{1}(\mathbb{N})$ is
module biflat by Theorem \ref{invflat}. Also we know for an
idempotent $e$ in an inverse semigroup $S$, $(e]=\{f\in E:
fe=ef=f\}$. For $S=\mathbb{N}$, for each $n\in \mathbb{N}$ we have
$(n]=\{m\in \mathbb{N}: m\geq n\}$. Hence $\mathbb{N}$ is not
uniformly locally finite (even not locally finite), so $\ell
^{1}(\mathbb{N})$ is not neither biprojective
 nor biflat (see the introduction).

 (ii) Let $\mathcal C$ be the bicyclic inverse
semigroup generated by $a$ and $b$, that is
$$\mathcal C=\{a^mb^n : m,n\geq 0 \},\hspace{0.2cm}(a^mb^n)^*=a^nb^m. $$
The set of idempotents of $\mathcal C$ is $E_{\mathcal C}=\{a^nb^n
: n=0,1,...\}$ which is totally ordered with the following order
$$a^nb^n \leq a^mb^m \Longleftrightarrow m \leq n.$$
It is shown in \cite{abe} that $\mathcal C/\approx$ is isomorphic to
the group of integers $\mathbb{Z}$, hence $\mathcal C$ is amenable.
Therefore $\ell ^{1}(\mathcal C)$ is module biflat, but not module
biprojective. It is easy to see $E_{\mathcal C}$ is not uniformly
locally finite, so $\ell ^{1}(\mathcal C)$ is neither biprojective
 nor biflat.

 (iii) Let $S$ be an amenable $E$-unitary
inverse semigroup with infinite number of idempotents (see
\cite{how} and \cite{pat}). Then ${\ell ^{1}}(S)$ is module
biflat.

(iv) If $S$ is a {\it Brandt semigroup} of an amenable group over
an infinite index set (see \cite{dun} and \cite{pat}), then ${\ell
^{1}}(S)$ is module biflat. \end{ex}

\end{document}